\documentclass[a4paper,12pt]{article}
\usepackage{amsmath,amssymb}
\def\al{\alpha}
\def\la{\lambda}
\def\dom{{\Omega^c}}
\def\domtx{{\omega(x,t)\,dx}}
\DeclareMathOperator\supp{supp}
\def\munu{M_1}
\def\mdoi{M_2}
\def\mtrei{M_3}
\def\munspe{M_4}
\def\mpatru{M_5}
\def\mcinci{M_6}
\def\mzece{M_7}
\def\msapte{M_8}
\def\mopt{M_9}
\def\mnoua{M_{10}}

\begin{document}
\title{Confinement of vorticity in two dimensional ideal incompressible exterior flow}
\author{D. Iftimie\\
M.C. Lopes Filho\thanks{Research supported in part by CNPq grant \#302.102/2004-3}  \\H.J. Nussenzveig Lopes\thanks{Research supported in part by CNPq grant \#302.214/2004-6}}
\date{}
\newtheorem{theorem}{Theorem}
\newtheorem{claim}[theorem]{Claim}
\newtheorem{lemma}[theorem]{Lemma}
\newtheorem{proposition}[theorem]{Proposition}
\newtheorem{definition}[theorem]{Definition}
\newtheorem{remark}[theorem]{Remark}
\newenvironment{proof}{\smallskip \noindent{\bf Proof}:}{      
               \hfill\rule{2mm}{3mm}\hspace{1in}\smallskip}
\newcommand{\R}{\mathbb{R}}
\newcommand{\C}{\mathbb{C}}
\newcommand{\loc}{{\scriptsize \mbox{loc}}}
\newcommand{\vare}{\varepsilon}
 
\maketitle
\begin{abstract}
In [Math. Meth. Appl. Sci. 19 (1996) 53-62], C. Marchioro examined the problem of vorticity confinement
in the exterior of a smooth bounded domain. The main result in Marchioro's paper is that solutions
of the incompressible 2D Euler equations with compactly supported nonnegative initial vorticity in the 
exterior of a connected bounded region have vorticity support with diameter growing at most like 
$\mathcal{O}(t^{(1/2)+\vare})$, for any $\vare>0$. In addition, if the domain is the exterior of a disk, 
then the vorticity support is contained in a disk of radius $\mathcal{O}(t^{1/3})$.
The purpose of the present article is to refine Marchioro's results. We will prove that, if the 
initial vorticity is even with respect to the origin, then the exponent for the exterior of the disk may 
be improved to $1/4$. For flows in the exterior of a smooth, connected, bounded domain we prove a confinement 
estimate with exponent $1/2$ (i.e. we remove the $\vare$) and in certain cases, depending on the harmonic part 
of the flow, we establish a logarithmic improvement over the exponent $1/2$. The main new ingredients in our 
approach are: (1) a detailed asymptotic description of solutions to the exterior Poisson problem near infinity, 
obtained by the use of Riemann mappings; (2) renormalized energy estimates and bounds on logarithmic moments 
of vorticity and  (3) a new {\it a priori} estimate on time derivatives of logarithmic perturbations of the moment 
of inertia.    
\end{abstract}

\section{Introduction}

Two-dimensional incompressible ideal flow can be described as the active transport of vorticity, see \cite{Constantin94}. 
Vorticity with changing sign may scatter through the divergent motion of soliton-like vortex pairs, see the discussion 
in \cite{ILN03b} and references there contained, but single-signed vorticity tends to rotate around, and spreads 
very slowly. Studying the rate at which single-signed vorticity spreads is the problem of vorticity confinement.

	In 1996, C. Marchioro presented some results concerning vorticity confinement in the case of exterior domain flow, see \cite{M}. Marchioro observed cubic-root confinement in the case of the exterior of a disk, i.e. single-signed, compactly supported, vorticity has its support contained in a space-time region
whose diameter grows like the cubic-root of time. This result follows from the proof of similar cubic-root confinement
obtained previously by Marchioro for full plane flow, see \cite{mar94}.  For flows in the exterior of a general connected domain, Marchioro proved $(1/2 + \vare)$-root confinement. The purpose of this article is to 
refine Marchioro's estimates. Our main result is unqualified square-root confinement for exterior flow. We improve this estimate to a logarithmic refinement of square-root confinement under certain conditions on the harmonic part of the flow. In addition, we prove almost fourth-root vorticity confinement in the exterior of a 
disk if the initial vorticity is even with respect to its center. Technically, we begin with the construction of
a conformal map between the exterior of a general domain and the exterior of the unit disk, which behaves nicely up to the boundary, taken from \cite{ILN03}. This conformal map is used to obtain explicit formulas for the Green's function of the exterior domain, the Biot-Savart law and the harmonic part of the velocity. We then prove {\it a priori} estimates: first, a renormalized energy estimate, next, in some cases, estimates on logarithmic moments of vorticity, and finally an estimate of linear growth in time for logarithmic perturbations of the moment of inertia. Finally, we use these {\it a priori} estimates to obtain our confinement results. 

      From a broad viewpoint, the problem of confinement is related with scaling.
Roughly speaking, scaling in an evolution equation is determined by the behavior in time of the radius of effective influence of a small localized perturbation. For a parabolic system, the scaling is $x \sim \sqrt{t}$ and for a hyperbolic system it is $x \sim t$. Incompressible ideal flow has interesting behavior at hyperbolic scaling, i.e. waves,  but this requires vortex pairs, and therefore vorticity changing sign, see \cite{ILN03b,ILN03c,ISG99} for details. One important  issue is whether there is a natural scaling associated with incompressible, ideal 2D flow with distinguished signed vorticity. Confinement estimates explore this issue of scaling, and therefore, are useful in studying the qualitative behavior of solutions. For example, confinement estimates have been used in the rigorous justification of point vortex dynamics as an asymptotic description of the dynamics of highly concentrated vorticity, see \cite{M98b,MP93} and in the results on vortex scattering in \cite{ILN03b,ILN03c}. In addition, the issue of confinement has attracted attention in
other contexts, such as confinement for slightly viscous flow, see \cite{M98a}, for axisymmetric flow, see \cite{BCM00,MM01,M99}, for the Vlasov-Poisson system, see \cite{CM00}, and for the quasigeostrophic system, see \cite{M96}.

The point of departure on vorticity confinement research is the 1994 article \cite{mar94},  by C. Marchioro, in which he 
proved cubic-root confinement for flows in the full plane. Marchioro used in an essential way the conservation of the 
moment of inertia, which is associated to the rotational symmetry of full plane flow.  This result was improved, independently by Ph. Serfati, \cite{Serfati98} and by D. Iftimie, T. Sideris and P. Gamblin, see \cite{ISG99}, to nearly fourth root in time. The improvement relied on using, in addition to the moment of inertia, conservation of the center of vorticity,
which is associated to the translational symmetry of the plane. However, scaling should be a robust qualitative property, and not dependent on the presence of symmetry. One of the main points of this article is to further explore the role that symmetry has in the problem of confinement. To this end we break symmetry by considering exterior domain flows, both in the exterior of a disk, where only translational symmetry is broken and in more general exterior domains, where translational and rotational symmetry are broken. We begin with the following two questions: (1) Is it possible to use the center of vorticity to improve Marchioro's estimate for the exterior domain in the same way that Serfati, Iftimie, Sideris and Gamblin did in the full-plane case? (2) Given that, without symmetry, there are no conserved quantities, is it still possible to find quantities which remain bounded in time and that could play the role of the moment of inertia for confinement estimates?  
Perhaps our most important result is negative -- we have not found coercive quantities like the moment of inertia for which
we could prove boundedness or slow growth. We proved boundedness for logarithmic moments of vorticity related to conservation
of energy and we proved that certain logarithmic perturbations of the moment of inertia grow at most linearly in time. Using these estimates we obtained small improvements over Marchioro's original estimate. Our results support the roughly parabolic natural scaling for incompressible ideal two-dimensional flow with distinguished signed vorticity, obtained originally
by Marchioro.        

The remainder of this paper is divided into five sections. In Section 2 we set up the problem and we collect
preliminary information on Riemann mappings and the Laplacian in an exterior domain. In Section 3 we derive
a renormalized energy estimate and we prove that, under appropriate hypothesis', boundedness of a 
logarithmic moment of vorticity can be deduced from the energy estimate. In Section 4 we obtain a linear-in-time 
estimate for logarithmic perturbations of quadratic moments of vorticity. In Section 5 we use the estimates
deduced in Sections 3 and 4 to obtain new confinement estimates for flows in the exterior
of a general connected domain. Finally, in Section 6 we consider flow induced by even vorticity in the exterior of a disk, obtaining nearly fourth-root confinement.

\section{Preliminaries on exterior domain flow}

	The purpose of this section is to collect information on exterior domain flow, particularly
the discussion on Riemann mappings and solutions of the Dirichlet and Poisson problems in an exterior
domain developed in \cite{ILN03}, supplementing the available results as required.

	We begin by discussing the asymptotic behavior, near infinity, of solutions to a 2D exterior domain Poisson equation. Let $\Gamma \subseteq \R^2$ be a smooth Jordan curve dividing the plane into a bounded connected 
component $\Omega$ and an unbounded connected component denoted $\Omega^c$. For $x \in \Gamma$, denote by $\hat{n}(x)$ the unit exterior normal to $\Omega^c$ at $x$. Note that both $\Omega$ and $\Omega^c$ are open. The Green's function for the Laplacian in $\Omega^c$ is denoted by
$G_{\Omega^c}$. We consider $\omega \in C^{\infty}_c (\Omega^c)$ and we introduce:
\begin{equation}
  \label{stream}
\psi(x) \equiv \int_{\Omega^c} G_{\Omega^c}(x,y) \omega(y) dy.  
\end{equation}
Let $S = \{|x| = 1\}$ and $D^c \equiv \{|x|>1\}$. Denote the inversion with respect to $S$ by 
$x \mapsto x^{\ast} \equiv x/|x|^2$. The Green's function in the case of the exterior of the unit disk can be written explicitly as:
\[ G_{D^c}(x,y) =  \frac{1}{2\pi}\log \frac{|x-y|}{|x-y^{\ast}||y|}. \]
It is easy to obtain detailed information on the asymptotic behavior of $\psi$ and its derivatives, in the case of the exterior of the disk, by means of the representation formula above. The objective of this section is to obtain similar 
information for general exterior domains. We begin with a version of Lemma 2.1 of \cite{ILN03}.

\begin{lemma} \label{confmap}
There exists a conformal mapping $T:\Omega^c \rightarrow D^c$ extending continuously up to the boundary, mapping $\Gamma$ to $S$. Furthermore, there exist $\beta$, a positive real number, and $h$, a bounded holomorphic function in $\Omega^c$,
such that 
\[T(z) = \beta z + h(z).\]
In addition, there exists a constant $\munu >0$ such that 
$|h^{\prime}(z)| \leq \munu /|z|^2$, $|h^{\prime\prime}(z)| \leq \munu /|z|^3$. Furthermore, if $T$ is regarded as a real mapping between $\Omega^c$ and $\{|z|>1\}$, then we have $\|DT\|_{L^{\infty}}\leq \munu $ and $\|DT^{-1}\|_{L^{\infty}} \leq \munu $. 
\end{lemma}

\begin{proof}
Most of the facts claimed above are either proved in Lemma 2.1 of \cite{ILN03} or in the remark following it.
One modification is that $\beta$ can be assumed positive, which we can obtain by composing $T$ with a rotation,
if necessary. The only other difference is the estimate on $h^{\prime\prime}$ which follows by differentiating
the relation $h^{\prime}(z) = -(1/z^2) g^{\prime}(1/z)$, together with the fact that, by construction, 
$g^{\prime}$ and $g^{\prime\prime}$ are bounded.

\end{proof}

	It is easy to verify that, if $x_0 \in D^c$ and $\phi$ satisfies $\Delta \phi = \delta(x-x_0)$ in a neighborhood of
$x_0$, then $\tilde{\phi} \equiv \phi \circ T$ satisfies $\Delta \tilde{\phi} = \delta(y-T^{-1}(x_0))$ in a neighborhood of $T^{-1}(x_0)$. We use this fact on the Green's function $G_{D^c}$ in order to write:
\begin{equation} \label{greenform}
G_{\Omega^c}(x,y) = \frac{1}{2\pi}\log \frac{|T(x)-T(y)|}{|T(x)-T(y)^{\ast}||T(y)|}. 
\end{equation}

We now formulate precisely the initial value problem for incompressible fluid flow in an exterior domain.  In this,
we follow Section 3.1 of \cite{ILN03}. Let us denote by $u = u(x,t) = (u_1(x_1,x_2,t),u_2(x_1,x_2,t))$ the velocity of an incompressible, ideal fluid flow in the exterior domain $\Omega^c$. We assume that $u$ is tangent to $\Gamma$ and $u \to 0$ when $|x| \to \infty$. The evolution of such a flow is governed by the Euler equations. We write the initial-boundary value problem as:
\begin{equation} \label{veleq}
\left\{ \begin{array}{ll}
u_t + u \cdot \nabla u = - \nabla p & \mbox{ in } \Omega^c \times (0,\infty)  \\
\mbox{div }u = 0 & \mbox{ in } \Omega^c \times [0,\infty) \\
u \cdot \hat{n} = 0 & \mbox{ in } \Gamma \times [0,\infty) \\
\lim_{|x| \to \infty}u = 0 & \mbox{ for } t \in [0,\infty)\\
u(x,0) = u_0(x) & \mbox{ in } \Omega^c,
\end{array} \right. \end{equation} 
where $p = p(x,t)$ is the scalar pressure. If $u_0$ is sufficiently smooth, global well-posedness of this problem was proved by K. Kikuchi in \cite{kikuchi83}. 

We introduce $\omega = \mbox{ curl }u$, the vorticity of the flow. Vorticity satisfies the transport equation:
\[ \omega_t + u \cdot \nabla \omega = 0, \mbox{ in } \Omega^c \times (0,\infty).\]
 
Our purpose is to reformulate problem (\ref{veleq}) in terms of vorticity.
In order to do this we require a version of the Biot-Savart law, which recovers velocity from vorticity. Since $\Omega^c$ is not simply connected, we must use Hodge-deRham theory (see \cite{warner71}). Recall that a {\it harmonic vector field} in $\Omega^c$ is a divergence-free, curl-free vector field tangent to $\Gamma$ and vanishing at infinity. By Hodge's Theorem, the vector space of harmonic vector fields, in our setting, is one-dimensional, see Section 2.3 in \cite{ILN03}. 
Therefore, every harmonic vector field is a multiple of a unique harmonic vector field, denoted $H_{\Omega^c}$, defined by requiring that  
\[ \int_{\Gamma} H_{\Omega^c} \cdot ds = 1.\]
(Throughout this paper circulation will be computed in the counterclockwise orientation.)
We recall identity (2.11) in \cite{ILN03}, which gives an explicit expression for 
$H_{\Omega^c}$ in terms of $T$. We have
\begin{equation} \label{harmform}
H_{\Omega^c}(x) = \frac{1}{2\pi} \nabla^{\perp} \log |T(x)| = \frac{1}{2\pi} 
\frac{(T(x)DT(x))^{\perp}}{|T(x)|^2} = \frac{1}{2\pi} 
\frac{DT^t(x)(T(x))^{\perp}}{|T(x)|^2}
\end{equation}

With this notation we can show that there exists $\alpha \in \R$ such that:
\begin{equation} \label{velchar}
u = \nabla^{\perp} \psi + \alpha H_{\Omega^c},
\end{equation}
where $\psi$ is the stream function introduced in \eqref{stream}.
Indeed, by \eqref{gradpsi}, $\nabla^{\perp} \psi$ vanishes at
infinity and, since $\psi$ vanishes on $\Gamma$, $\nabla^{\perp}\psi$
is tangent to $\Gamma$. Clearly $\nabla^{\perp}\psi$ is
divergence-free and its curl is $\omega$. Thus $u-\nabla^{\perp}\psi$
is a harmonic vector field, which must then be a real multiple of
$H_{\Omega^c}$, (see \cite[Proposition 2.1]{ILN03}). 

In the language of Hodge theory, the vector field $\alpha H_{\Omega^c}$ is called the {\it harmonic part} of the flow $u$. In principle this harmonic part is time-dependent, but as a consequence of Kelvin's Circulation Theorem, it is actually a constant of motion. 

\begin{lemma} \label{alphaconst} 
If $u$ is a solution of (\ref{veleq}) then $\alpha$ is  constant in time. 
\end{lemma}

This is Lemma 3.1 in \cite{ILN03}. More precisely, 
let $u_0$ be such that $\mbox{ curl }u_0$ is compactly supported and
set $\omega_0 \equiv \mbox{ curl }u_0$. Then 
\begin{equation}
  \label{alpha}
\alpha \equiv  \int_{\Gamma} u_0 \cdot ds + \int_{\Omega^c} \omega_0 dx.  
\end{equation}

With this notation the vorticity formulation of the initial-boundary value problem (\ref{veleq}) is: 
\begin{equation} \label{ivpvorteq}
\left\{ \begin{array}{ll}
\omega_t + u \cdot \nabla \omega  =  0 & \mbox{ in } \Omega^c \times (0,\infty)  \\
u = \nabla^{\perp} \psi + \alpha H_{\Omega^c} & \mbox{ in } \Omega^c \times [0,\infty) \\
\psi(x,t) = \int_{\Omega^c} G_{\Omega^c}(x,y) \omega(y,t) dy & \mbox{ in } 
\Omega^c \times [0,\infty)\\
\omega(x,0) = \omega_0(x) & \mbox{ in } \Omega^c.
\end{array} \right. 
\end{equation}
 
The fact that vorticity is transported by a divergence free vector field implies that
its $L^p$-norm is conserved for any $1 \leq p \leq \infty$.

We will also require information on the time-dependent stream function $\psi=\psi(\cdot,t)$, 
defined by \eqref{stream} with $\omega=\omega(\cdot,t)$. 

\begin{lemma} \label{bddstream} There exists a constant $C=C(t)>0$, depending on the diameter of the
support of $\omega(\cdot,t)$ such that
\[\sup_{x\in\Omega^c} |\psi(x,t)| \leq C(t).\]
\end{lemma}

\begin{proof}
We repeat the argument leading to relation (4.10) in \cite{ILN03}, substituting $\omega_t$
by $\omega$, to obtain
\[ \left| \psi(x,t) + \frac{1}{2\pi} \int_{\Omega^c} \log|T(y)| \omega(y,t) \, dy \right| = \mathcal{O}
\left(\frac{1}{|x|}\right). \]
Since $\omega(\cdot,t)$ has compact support, the conclusion follows.
\end{proof} 

In addition, we recall estimates (2.8) and (4.11) of \cite{ILN03}:
\begin{equation}
  \label{gradpsi}
|\nabla\psi(x,t)| \leq \frac{C(t)}{|x|^2}\qquad\text{and}\qquad  |\nabla\psi_t(x,t)| \leq \frac{C(t)}{|x|^2},
\end{equation}
where, again, $C(t)>0$ depends on the diameter of the support of $\omega(\cdot,t)$.

\section{Generalized energy and logarithmic moment}

	The purpose of this section is to derive a new {\it a priori} bound on the 
logarithmic moment of vorticity, in terms of global conserved quantities of the flow. 
 Our point of departure is the exact conservation of an energy-like quantity which we 
will call {\it generalized energy}. To define this quantity, let us consider $u = u(x,t)$ 
a smooth solution of problem (\ref{veleq}) with compactly supported vorticity. 
Using \eqref{gradpsi} and \eqref{velchar} we conclude that:
\begin{equation} \label{vdef}
v(x,t) \equiv u(x,t) - \alpha H_{\Omega^c}(x) = (\nabla^{\perp} \psi)(x,t) = {\cal O}(1/|x|^2), \mbox{ when } |x| \to \infty. 
\end{equation}

\begin{definition} \label{genen}
We define the generalized energy $E$ by:
\[ E \equiv \int_{\Omega^c} (|v|^2 + 2 \alpha H_{\Omega^c} \cdot v) dx. \]
\end{definition}

The observation on the asymptotic behavior of $v$ as $|x| \to \infty$ allows us to
conclude that $E$ is finite. Indeed $v$ is $O(|x|^{-2})$ and from
\eqref{harmform} and Lemma \ref{confmap} the vector field
$H_{\Omega^c} $ is $O(|x|^{-1})$ so that the integrand above is
$O(|x|^{-3})$ which is an integrable function at infinity.

\begin{proposition} \label{consgenen}
The generalized energy $E$ is a constant of motion for smooth flows on $\Omega^c$.
\end{proposition}

\begin{proof}
We begin this proof by noting that, thanks to \eqref{gradpsi}, for each fixed time $t \geq 0$ we have  
\[ |u_t(x,t)| = |v_t(x,t)| = |\nabla^\perp\psi_t(x,t)| = {\cal O}(1/|x|^2), \mbox{ as } |x| \to \infty. \]

The construction of $H_{\Omega^c}$ given in \eqref{harmform}  implies that $|DH_{\Omega^c}| = {\cal O}(1/|x|^2)$ as $|x| \to \infty$.
Hence, the same argument used above for $u_t$ yields the conclusion that 
$|Du|(x,t) = {\cal O}(1/|x|^2)$ as well. Equation (\ref{veleq})
implies that the same conclusion
holds for the behavior of $|\nabla p|$ at infinity. By integrating along rays, recalling that $p$ is only defined up to a constant, we can
conclude that $p = {\cal O}(1/|x|)$ as $|x| \to \infty$.

We use the velocity formulation (\ref{veleq}), and Lemma \ref{alphaconst} 
to compute the time derivative of $E$ as follows:
\begin{align*}
 \frac{dE}{dt} 
&= \int_{\Omega^c} (2 v \cdot v_t + 2\alpha H_{\Omega^c}\cdot v_t) dx\\
&= 2\lim_{R\to \infty}
 \int_{B(0;R) \setminus \Omega} \bigl[-v(u\cdot\nabla u + \nabla p) - \alpha H_{\Omega^c}(u\cdot\nabla u + \nabla p)\bigr]dx,\\
&= -2\lim_{R\to \infty} \left( \int_{B(0;R)\setminus \Omega} u\cdot(u\cdot \nabla u) dx + 
\int_{B(0;R) \setminus \Omega}u \cdot \nabla p dx \right) \\
&= -2\lim_{R\to \infty}({\cal I}_1 + {\cal I}_2). 
\end{align*}
 
We estimate these integrals:
\begin{gather*}
{\cal I}_1 = \int_{B(0;R) \setminus \Omega} \mbox{ div }
\frac{u|u|^2}{2} dx = \int_{|x| = R}\frac{|u|^2}{2}u\cdot\frac{x}{R}dS
+ \int_{\Gamma}\frac{|u|^2}{2}u\cdot\hat{n}dS = {\cal O}\left(
  \frac{1}{R^2} \right),\\
\intertext{and}
{\cal I}_2 = \int_{B(0;R) \setminus \Omega} \mbox{ div } (up) dx =
 \int_{|x| = R}p (u\cdot\frac{x}{R})dS  + \int_{\Gamma}p(u\cdot\hat{n})dS = {\cal O}\left(\frac{1}{R}\right).  
\end{gather*}
Thus, $dE/dt = 0$, as we desired.
\end{proof}

Instead of estimating the logarithmic moment of vorticity directly we will estimate
a quantity that resembles the logarithmic moment, but which is adapted
to the geometry of the domain under consideration. Recall that $T$ is the conformal map that takes $\Omega^c$ into the exterior of the unit disk. We introduce the {\it modified logarithmic moment} of vorticity by:
\begin{equation} \label{logmomm}
L(t) \equiv \frac{1}{2\pi}\int_{\Omega^c} \left( \log|T(y)| \right)\omega(y,t) dy.
\end{equation} 

Throughout the remainder of this paper we will assume that $\omega_0$ is smooth, nonnegative and compactly supported. 
Since $\omega$ is a solution of a transport equation with smooth velocity, it remains compactly supported and 
nonnegative for positive time. The total mass of vorticity is a conserved quantity and we denote it by
\[ m \equiv \int_{\Omega^c}\omega_0(x) dx = \int_{\Omega^c} \omega(x,t) dx. \]

\begin{theorem} \label{logmom}
If either $\alpha < 0$ or $\alpha > m$ then there exists a constant $\mdoi >0$ such that $L(t) \leq \mdoi .$
\end{theorem}

\begin{proof}
We first rewrite the generalized energy in terms of vorticity. By \eqref{vdef} and \eqref{harmform} 
we have
\begin{align*}
E = \int_{\Omega^c} |v|^2 + 2 \alpha H_{\Omega^c} \cdot v dx
&= \int_{\Omega^c} \nabla^{\perp} \psi \cdot v dx + \frac{\alpha}{\pi} 
\int_{\Omega^c} (\nabla^{\perp} \log |T(x)|) \cdot v dx \\
&= -\int_{\Omega^c} \psi \omega dx - \frac{\alpha}{\pi} 
\int_{\Omega^c} \log |T(x)| \omega dx,
\end{align*}
where the boundary terms vanish since, on $\Gamma$, $\psi$ and $\log|T|$ vanish, and 
at infinity we have, by Lemma \ref{bddstream}, that $\psi$ is bounded and, by \eqref{vdef}, that $v$ decays
like $O(|x|^{-2})$.

Using \eqref{stream} and \eqref{greenform}, we rewrite the energy in the following way:
\[-2\pi E = \iint \log \frac{|T(x)-T(y)|}{|T(x)-T(y)^{\ast}||T(y)|} 
\omega(x,t)\omega(y,t)\,dxdy + 2\alpha \int \log |T(x)|\omega(x,t) dx, \]
where, in the three integrals, the domain of integration is $\Omega^c$. Therefore we have
\[
-2\pi E = \iint \log |T(x)-T(y)| \omega(x,t)\omega(y,t) dxdy + 2(\alpha-m) \int \log |T(x)| \omega(x,t)dx +\]
\begin{equation} \label{eqq1}
+ \iint \log \frac{|T(x)|}{|T(x) - T(y)^{\ast}|} \omega(x,t)\omega(y,t)dxdy.
\end{equation}

We begin by observing that the last integral on the right-hand-side of \eqref{eqq1} is bounded independently of 
time. Indeed, if $z_1,z_2 \in \R^2$ with $|z_2|>1$ and $|z_1|>2$ then we have that 
\[2|z_1| \geq |z_1| + |z_2^{\ast}| \geq |z_1 - z_2^{\ast}| \geq |z_1| - |z_2^{\ast}| \geq |z_1|/2.\]
Applying this inequality with $z_1 = T(x)$ and $z_2 = T(y)$ yields
\begin{align*}
\biggl| \iint  \log 
&\frac{|T(x)|}{|T(x) - T(y)^{\ast}|}
  \omega(x,t)\omega(y,t)dxdy \biggr| \\
&\leq  m^2 \log 2 + 
\left| \int_{1 \leq |T(x)| \leq 2} \int \log \frac{|T(x)|}{|T(x) - T(y)^{\ast}|} \omega(x,t)\omega(y,t)dxdy
\right| \\
&\leq 2 m^2 \log 2 + 
\int_{1 \leq |T(x)| \leq 2} \int \left|\log |T(x) -
  T(y)^{\ast}|\right| \omega(x,t)\omega(y,t)dxdy\\
&\leq 2 m^2 \log 2 + \|\omega\|_{L^{\infty}} \int \left(
\int_{1 \leq |T(x)| \leq 2}  \left|\log |T(x) - T(y)^{\ast}|\right|
dx\right) \omega(y,t)dy\\
&\leq 2 m^2 \log 2 +  m \|\omega\|_{L^{\infty}} 
(\sup_{1 \leq |\eta| < 2} |\det DT(T^{-1}(\eta))|^{-1}) \int_{|\eta|
  \leq 3} |\log |\eta||d\eta,
\end{align*}
where we used the boundedness of derivatives of $T$ and $T^{-1}$
stated in Lemma \ref{confmap}. 
This establishes the desired time-independent bound. 

Next we decompose $\Omega^c \times \Omega^c$ as
\[\Omega^c \times \Omega^c = (\Omega^c \times \Omega^c \cap \{|T(x)-T(y)| \leq 1\}) \cup
 (\Omega^c \times \Omega^c \cap \{|T(x)-T(y)| > 1\}) \equiv A_1 \cup A_2.\]

An argument similar to the one carried out above
implies that 
\begin{multline*}
 \left| \iint_{A_1} \log |T(x)-T(y)| \omega(x,t) \omega(y,t) dx dy
 \right| \\
\leq m \|\omega\|_{L^{\infty}}
\sup_{|\eta|>1}|\det DT(T^{-1}(\eta))|^{-1} \int_{|\eta|\leq 1} |\log(|\eta|)| d\eta.
\end{multline*}

These estimates, together with \eqref{eqq1} and conservation of energy imply that
\begin{equation}
  \label{it}
\mathcal{I}(t) \equiv \iint_{A_2} \log |T(x)-T(y)| \omega(x,t) \omega(y,t) dx dy + 2(\alpha - m) \int \log|T(x)| \omega(x,t) dx  
\end{equation}
is bounded independently of time.
Now, if $\alpha > m$ both terms are nonnegative and the logarithmic moment bound follows immediately.

Lastly, we treat the case $\alpha < 0$. Since $|T(x)|$, $|T(y)| \geq 1$ it follows that  
\[ \log|T(x) - T(y)| \leq \log|T(x)| + \log|T(y)| + \log 2. \]
Therefore, 
\[\left| \iint_{A_2} \log |T(x)-T(y)| \omega(x,t) \omega(y,t) dx dy \right| \leq m^2\log 2  +   
2m \int \log|T(x)| \omega(x,t) dx,\]
which, together with the boundedness of $\mathcal{I}$ concludes the proof.
\end{proof}

We conclude this section with an estimate which applies to the extreme case $\alpha = 0$.
We do not have a logarithmic moment bound in this case, but we can prove another estimate
which we will be able to use in place of a logarithmic moment bound in the analysis
that follows. 

\begin{lemma} \label{alphazero}
Assume that $\alpha = 0$. Then there exists a constant $\mtrei >0$ such that
\[ \iint_{\Omega^c\times\Omega^c} \log \left( \frac{|T(x)||T(y)|}{|T(x) - T(y)|} \right) \omega(x,t)
\omega(y,t) dxdy \leq \mtrei . \]
\end{lemma}

\begin{proof}
First recall the definition of $\mathcal{I}(t)$ given in \eqref{it}. Next observe that, if $\alpha = 0$, then the 
integral we wish to estimate is equal to 
\begin{equation*}
-\mathcal{I}(t) -\iint_{A_1} \log |T(x)-T(y)| \omega(x,t) \omega(y,t) dx dy,  
\end{equation*}
in this case. Finally, the proof of Theorem \ref{logmom} shows that both of these
terms are bounded independently of the choice of $\alpha$.
\end{proof}

\begin{remark}
If either $\al\leq0$ or $\al>m$, then there exists a constant $\munspe > 0$ such that
\begin{equation}
  \label{logfin}
  \iint_{\Omega^c\times\Omega^c} \left[ \log \min (|T(x)|,|T(y)|) \right] \omega(x,t)
\omega(y,t) dxdy \leq \munspe.
\end{equation}
\end{remark}
Indeed, if  $\al<0$ or $\al>m$ then the relation above is an immediate
consequence of Theorem \ref{logmom}. If $\al=0$, then \eqref{logfin}
follows from Lemma \ref{alphazero} after we observe that
\begin{equation*}
 \log\left(\frac{|T(x)||T(y)|}{|T(x) - T(y)|} \right) \geq \log \left( \frac{\min (|T(x)|,|T(y)|)}{2} \right). 
\end{equation*}

\section{Other moment bounds}
\label{secmom}

The purpose of this section is to derive bounds on logarithmic perturbations of the
moment of inertia, which will later be used to improve Marchioro's confinement estimate. 
We begin with a technical lemma regarding the conformal map $T$.

\begin{lemma} \label{loops}
Let $x,y \in \Omega^c$. We have
\begin{equation} \label{loops1}
|T(x) \cdot (T(y))^{\perp}| \leq \min\bigl\{|T(x)|,|T(y)|\bigr\}|T(x) - T(y)|,
\end{equation}
and  there exists $\mpatru >0$ such that
\begin{equation} \label{loops2}
||T^{\prime}(x)|^2 - |T^{\prime}(y)|^2| \leq \mpatru \frac{|T(x) - T(y)|}{\min\{|T(x)|,|T(y)|\}^2}
\end{equation}
\end{lemma}

\begin{proof} Let us first prove \eqref{loops1}. We observe that
\[|T(x) \cdot (T(y))^{\perp}| = |(T(x) - T(y)) \cdot (T(y))^{\perp}| \leq |T(x) - T(y)||T(y)|. \]
Similarly,
\[|T(x) \cdot (T(y))^{\perp}| = |T(x)\cdot (T(x) - T(y))^{\perp}| \leq |T(x) - T(y)||T(x)|. \]
Hence, \eqref{loops1} follows from both these inequalities.

Next we prove \eqref{loops2}.  We use Lemma \ref{confmap} to compute
\begin{align*}
||T^{\prime}(x)|^2 - |T^{\prime}(y)|^2| 
&= |2 \beta \mbox{Re} (h^{\prime}(x) - h^{\prime}(y)) + |h^{\prime}(x)|^2 - 
|h^{\prime}(y)|^2| \\
&\leq C |h^{\prime}(x) - h^{\prime}(y)| + ||h^{\prime}(x)|^2 - |h^{\prime}(y)|^2|\\
&\leq C|h^{\prime}(x) - h^{\prime}(y)|(1 + |h^{\prime}(x)| +
|h^{\prime}(y)|) \\
&\leq C|h^{\prime}(x) - h^{\prime}(y)|,
\end{align*}
since, by Lemma \ref{confmap}, $|h^{\prime}|$ is bounded.

Therefore, to conclude our proof it is enough to show that $|h^{\prime}(x) - h^{\prime}(y)|$ is bounded 
by the right-hand-side of \eqref{loops2}.

We first remark that a uniform bound on the first derivatives of a map defined on $\Omega^c$ implies a global
Lipschitz bound. Next note that, since, by Lemma \ref{confmap}, $h^{\prime\prime}$ is bounded, the argument
above implies that $h^{\prime}$ is globally 
Lipschitz and therefore $|h^{\prime}(x) - h^{\prime}(y)| \leq C |x-y|$. Furthermore, 
$|x-y| \leq C |T(x)-T(y)|$, with $C$ being the Lipschitz constant associated to $T^{-1}$, and using again the global
bound on the derivatives of $T^{-1}$. Therefore, if either $x$ or $y$ is contained in a ball of radius $R$, then 
\[ \frac{|h^{\prime}(x) - h^{\prime}(y)|}{|T(x) - T(y)|} \leq C \leq C(R) \frac{1}{\min\{|T(x)|,|T(y)|\}^2}.\]  
 It remains only to prove this inequality assuming both $|x| > R$ and $|y| > R$. In this case, there exists a smooth path $\gamma:[0,1] \to \{|z| \geq R\}$, whose length is less than $2|x-y|$, and such that 
$\gamma(0) = x$, $\gamma(1) = y$. Suppose first that $|x| \leq |y|$.
We can also assume that $\gamma(s) \in \{z \;|\; |x| \leq |z| \leq |y|\}$ for all $s \in [0,1]$. We then have:
\begin{multline*}
|h^{\prime}(y) - h^{\prime}(x)| = 
\left| \int_0^1 \frac{d}{dt} h^{\prime}(\gamma(t)) dt \right| 
\leq \int_0^1 |h^{\prime\prime}(\gamma(t))||\gamma^{\prime}(t)| dt \\
\leq \frac{C}{|x|^2} \mbox{ Length}(\gamma)  \leq \frac{C|x-y|}{(\min{|x|,|y|})^2}
\leq \frac{C|T(x)-T(y)|}{(\min{|x|,|y|})^2},  
\end{multline*}
where we have used, in the second inequality, the decay estimate for $|h^{\prime\prime}|$ from Lemma \ref{confmap}.

A similar argument holds for $|y| \leq |x|$. This completes the proof.
\end{proof}

The final result in this section is an estimate on the growth of logarithmic perturbations of the
moment of inertia.

\begin{theorem} \label{keythm}
There exists a constant $\mcinci >0$ such that
\[\int_{\Omega^c} |T(x)|^2 \left( \log|T(x)| \right) \omega(x,t) dx \leq \mcinci (1+t),\]
for all $t \geq 0$. Moreover, if either $\alpha \leq 0$ or $\alpha > m$ then 
\[\int_{\Omega^c} |T(x)|^2 \left( \log^2|T(x)| \right) \omega(x,t) dx \leq \mcinci (1+t),\]
for all $t\geq 0$.
\end{theorem}

\begin{proof}
Let $\sigma : (1,\infty) \to \R$ be a smooth function and define
\[\mathcal{J}_{\sigma}(t) \equiv \int_{\Omega^c} \sigma(|T(x)|^2)\omega(x,t) dx.\]
Clearly it is enough to prove that, for the appropriate choice of $\sigma$
(either $\sigma(s)=s\log s$ or $\sigma(s)=s\log^2 s$), the time-derivative
of $\mathcal{J}_{\sigma}$ is bounded independently of time. We estimate directly
\begin{align*}
\mathcal{J}_{\sigma}^{\prime}(t) 
&= -\int_{\Omega^c} \sigma(|T(x)|^2) \mbox{ div}( u\omega)(x,t)dx \\
&=\int_{\Omega^c} \nabla(\sigma(|T(x)|^2)) \cdot u(x,t) \omega(x,t) dx \\
&= \int_{\Omega^c} \nabla(\sigma(|T(x)|^2)) \cdot v(x,t) \omega(x,t) dx + \alpha \int_{\Omega^c} \nabla(\sigma(|T(x)|^2)) \cdot 
H_{\Omega^c}(x) \omega(x,t) dx\\
&\equiv I_1 + I_2.
\end{align*}
We begin by observing that $I_2=0$. Indeed,
\[H_{\Omega^c}(x) = \frac{1}{2\pi}\nabla^{\perp}\log(|T(x)|) = \frac{1}{4\pi} \nabla^{\perp} \log |T(x)|^2.\]
Clearly, $\sigma(|(T(x)|^2)$ and $\log (|T(x)|^2)$ are functionally dependent, which implies that their gradients
are proportional everywhere. Therefore, the integrand in $I_2$ vanishes identically.

Next, we estimate $|I_1|$. We use the explicit expression of the Biot-Savart kernel in \cite{ILN03}, equation (2.5)
to write
\[v(x,t) = \int_{\Omega^c} \left(\frac{[(T(x)-T(y))DT(x)]^{\perp}}{2\pi|T(x)-T(y)|^2}
- \frac{[(T(x)-T(y)^{\ast})DT(x)]^{\perp}}{2\pi|T(x)-T(y)^{\ast}|^2}\right) \omega(y,t) dy.\]

We write $I_1 \equiv I_{11} - I_{12}$ with

\[I_{11} \equiv \iint_{\Omega^c \times \Omega^c} \nabla(\sigma(|T(x)|^2)) \cdot 
\frac{[(T(x)-T(y))DT(x)]^{\perp}}{2\pi|T(x)-T(y)|^2} \omega(x,t)\omega(y,t) dx dy, \]
and 
\[I_{12} \equiv \iint_{\Omega^c \times \Omega^c} \nabla(\sigma(|T(x)|^2)) \cdot 
\frac{[(T(x)-T(y)^{\ast})DT(x)]^{\perp}}{2\pi|T(x)-T(y)^{\ast}|^2} \omega(x,t)\omega(y,t) dx dy, \]

We estimate $I_{12}$ first.  We further decompose
\[I_{12} = \int_{1 \leq |T(x)| \leq 2} \int_{\Omega^c} \cdots dydx + \int_{2 \leq |T(x)|} \int_{\Omega^c} \cdots dydx
\equiv I_{121} + I_{122}. \]

We have
\begin{multline*}
|I_{121}| \leq C m \sup_{|\eta| \leq 1}\left( \int_{1 \leq |T(x)| \leq
    2} \frac{\omega(x,t)}{|T(x) - \eta|} dx\right)
\leq C m \|\omega\|_{L^{\infty}} \sup_{|\eta|\leq 1} \left(\int_{1
    \leq |\zeta| \leq 2} 
\frac{d\zeta}{|\zeta - \eta|}\right)\\
\leq C m \|\omega\|_{L^{\infty}}\int_{|\zeta| \leq 3} \frac{d\zeta}{|\zeta|} \leq C m \|\omega\|_{L^{\infty}}.  
\end{multline*}

Next we treat $I_{122}$. Suppose first that $\sigma(s) = s \log^2 s$. Then
\begin{align*}  
\nabla(\sigma(|T(x)|^2)) 
&= 2\sigma^{\prime}(|T(x)|^2) T(x) DT(x) \\ 
&= 4[\log^2(|T(x)|) + \log(|T(x)|)]T(x)DT(x) \equiv [o(|T(x)|)]T(x)DT(x).
\end{align*}
Hence, using that $T(x)=\beta x + \mathcal{O}(1)$ and $DT(x) = \beta Id + \mathcal{O}(1/|x|^2)$ (see Lemma 
\ref{confmap}), we obtain, 
when $|T(x)|>2$, that:
\[|I_{122}| = \left|C\int_{|T(x)|\geq 2}\int_{\Omega^c} [o(|T(x)|)]T(x)DT(x) \cdot 
\frac{[(T(y)^{\ast})DT(x)]^{\perp}}{2\pi|T(x)-T(y)^{\ast}|^2} \omega(x,t)\omega(y,t) dx dy\right| \]
\[\leq m^2 \sup_{\{|T(x)|\geq 2\}} \frac{o (|T(x)|)}{|T(x)|},\]
which shows that $|I_{122}|$ is bounded independently of time. 

A similar argument may be used in case  $\sigma(s) = s \log s$.

Finally, we must estimate $I_{11}$. First we observe that if we consider the holomorphic map $T$ as a 
real map from $\R^2$ to itself, we have that $(v DT)^{\perp} =
v^{\perp} DT^t$. We also note that, by the Cauchy-Riemann relations, the matrix $DT$ has the form
\[\left[ \begin{array}{rr} a & b \\ -b & a \end{array} \right],\]
and therefore, 
$DT(x) DT^t(x) = \det(DT(x)) Id = |T^{\prime}(x)|^2 Id$. Consequently,
\begin{align*}
 \nabla(\sigma(|T(x)|^2)) \cdot [(T(x)-&T(y))DT(x)]^{\perp} \\
&=2\sigma^{\prime}(|T(x)|^2) [T(x) DT(x)]\cdot [(T(x)-T(y))^\perp DT(x)]\\
&=2\sigma^{\prime}(|T(x)|^2) [T(x) DT(x)DT^t(x)]\cdot (T(x)-T(y))^\perp\\
&=-2\sigma^{\prime}(|T(x)|^2) |T^{\prime}(x)|^2 T(x) \cdot T(y)^\perp
\end{align*}
Plugging this expression in the definition of $I_{11}$ we get
\begin{align*}
I_{11}
&=  -\iint \sigma^{\prime}(|T(x)|^2)
|T^{\prime}(x)|^2 \frac{T(x) \cdot T(y)^\perp}{\pi|T(x)-T(y)|^2}
\omega(x,t)\omega(y,t) dx dy\\
&= \iint \bigl[ \sigma^{\prime}(|T(y)|^2)
|T^{\prime}(y)|^2 -\sigma^{\prime}(|T(x)|^2)
|T^{\prime}(x)|^2\bigr] \frac{T(x) \cdot T(y)^\perp}{2\pi|T(x)-T(y)|^2}
\omega(x,t)\omega(y,t) dx dy\\
&\equiv \iint A(x,y) \omega(x,t)\omega(y,t) dx dy.
\end{align*}
To obtain the second line in the previous relation, we observed 
that  
\[T(x) \cdot T(y)^\perp=-T(y) \cdot T(x)^\perp \]
and we symmetrized, by adding half of the integrand with half of the same expression but with $x$ and $y$
interchanged. 

We see that, to conclude the proof, it is
sufficient to prove that if $\sigma'(s)=1+\log s$, then $|A|$ is bounded by
a constant and, in view of \eqref{logfin},
if $\sigma'(s)=\log^2s+2\log s$, then $|A|$ is bounded by $C+C\log\min(|T(x)|,|T(y)|)$.
Observe next that, since $A(x,y)$ is symmetric in $x$ and $y$, we can
assume for example that $|T(y)|\leq |T(x)|$. Using Lemma \ref{loops} we first bound
\begin{align*}
  2\pi |A(x,y)|
&= \biggl| \Bigl\{\bigl[ \sigma^{\prime}(|T(x)|^2)- \sigma^{\prime}(|T(y)|^2)\bigr]
|T^{\prime}(x)|^2 \\
&\hskip 4cm +\sigma^{\prime}(|T(y)|^2)
\bigr[|T^{\prime}(x)|^2-|T^{\prime}(y)|^2\bigr] \Bigr\}\frac{T(x)
  \cdot T(y)^\perp}{|T(x)-T(y)|^2}\biggr|\\
&\leq \bigl| \sigma^{\prime}(|T(x)|^2)-
\sigma^{\prime}(|T(y)|^2)\bigr|\|DT\|_{L^\infty}^2   \frac{|T(y)|}{|T(x)-T(y)|} +\sigma^{\prime}(|T(y)|^2)\frac{\mpatru}{|T(y)|}\\
&\leq C \bigl| \sigma^{\prime}(|T(x)|^2)-
\sigma^{\prime}(|T(y)|^2)\bigr|
\frac{|T(y)|}{|T(x)-T(y)|} +C,
\end{align*}
where the last inequality holds since, for both choices  $\sigma'(s)=1+\log s$ and $\sigma'(s)=\log^2s+2\log s$, 
we have $|\sigma^{\prime}(T(y))| = o(|T(y)|)$.

\bigskip

Let us now consider the case $\sigma'(s)=1+\log s$. We will use the elementary inequality
\begin{equation} \label{elemineq1}
|\log a-\log b|\leq \frac{|a-b|}{\min(a,b)}, \;\;\; a, \, b>0,
\end{equation} 
to observe that
\begin{equation*}
\bigl| \sigma^{\prime}(|T(x)|^2)-
\sigma^{\prime}(|T(y)|^2)\bigr|
\frac{|T(y)|}{|T(x)-T(y)|}\leq 2  
\end{equation*}
which implies that $A$ is bounded. 

Finally, if  $\sigma'(s)=\log^2s+2\log
s$, we use the inequality
\begin{equation} \label{elemineq2}
|\log^2 a+2\log a-\log^2 b-2\log b|\leq2|a-b|\frac{[1+\log\min (a,b)]}{\min (a,b)}
\end{equation}
together with \eqref{elemineq1} to deduce that
\begin{equation*}
\bigl| \sigma^{\prime}(|T(x)|^2)-
\sigma^{\prime}(|T(y)|^2)\bigr|
\frac{|T(y)|}{|T(x)-T(y)|}\leq 12+ 8 \log|T(y)|  \leq  C [1+\log \min(|T(x)|,|T(y)|)].  
\end{equation*}

The proof is completed.
\end{proof}

\section{Dispersion of Vorticity}

In this section, we use the bounds on the logarithmic moments of
inertia obtained in Section \ref{secmom} to deduce some confinements results for
the support of vorticity. We start with a useful estimate.
 
\begin{lemma}\label{velbd}
Let $S\subset\Omega^c$ and $\zeta:S\rightarrow{\mathbb{R}}^+$ be a function belonging
to $L^1(S)\cap L^\infty(S)$. There exists a constant $C>0$ such that
\begin{gather*}
\int_S\frac{\zeta(y)}{|x-T(y)|}\, dy
\leq C\|{\zeta}\|_{L^{1}(S)}^{1/2}\|{\zeta}\|_{L^\infty(S)}^{1/2}\qquad
\forall x\in\R^2\\
\intertext{and}
\int_S\frac{\zeta(y)}{|x-T(y)^*|}\, dy
\leq
C\|{\zeta}\|_{L^{1}(S)}^{1/2}\|{\zeta}\|_{L^\infty(S)}^{1/2}+C\|{\zeta}\|_{L^1(S)}\qquad
\forall |x|>1.
\end{gather*}
\end{lemma}
\begin{proof}
The case $T=Id$ of the first relation is  \cite[Lemma 2.1]{ISG99}. The general case follows
from this particular case in the following manner: 
\begin{align*}
 \int_S\frac{\zeta(y)}{|x-T(y)|}\, dy
&\leq C\int_{T(S)} \frac{1}{|x-z|}
\zeta(T^{-1}(z)) dz\\
&\leq C\|\zeta\|_{L^\infty(S)}^{\frac12}\Bigl(\int_{T(S)} 
\zeta(T^{-1}(z)) dz\Bigr)^{\frac12}\\
&\leq C\|\zeta\|_{L^\infty(S)}^{\frac12}
\Bigl(\int_{S} 
\zeta(y) dy\Bigr)^{\frac12}.
\end{align*}
We used above that the Jacobians of $T$ and $T^{-1}$ lie between two
positive constants to deduce the first and third lines; we also
used the particular case $T=Id$ to obtain the second line.

The second inequality follows from the first one after we observe that if
$|x|\geq2$ or $|T(y)|\geq2$ then $|x-T(y)^*|\geq\frac12$, while if
$|x|\leq2$ and $|T(y)|\leq2$ then $|x-T(y)^*|=|x| |x^*-T(y)|/|T(y)|\geq |x^*-T(y)|/2$.
\end{proof}

We prove now the following theorem.

\begin{theorem}\label{th1}
Assume that the initial vorticity $\omega_0$ is nonnegative, bounded,
compactly supported and let
$\omega(x,t)$ be the solution of \eqref{ivpvorteq}. There exists a
constant $\mzece >0$ such that the support of $\omega(x,t)$ is confined in
the set $|x|\leq \mzece  (1+t)^{\frac12}$. Furthermore, if either $\al\leq0$ or
$\al>m$ then  the support of $\omega(x,t)$ is confined in
the set $|x|\leq \mzece  (1+t)^{\frac12} [\log(2+ t)]^{-\frac{1}{4}}$. 
\end{theorem}
\begin{proof}
We introduce the following parameter $\theta$: we set $\theta=2$ if either 
$\al\leq0$ or $\al>m$ and $\theta=1$ otherwise. According to Theorem
\ref{keythm}, there exists a constant $\mcinci >0$ such that   
  \begin{equation}
    \label{momint}
    \int_\dom |T(x)|^2\log^\theta|T(x)|\domtx\leq \mcinci(1+t)
  \end{equation}
 for all $t$. The conclusion of Theorem \ref{th1} is tantamount to
 proving that there exists a constant $\mzece >0$ such
 that the support of $\omega(x,t)$ is confined in
the set $|x|\leq \mzece  (1+t)^{\frac12} [\log(2+ t)]^{\frac{1-\theta}4}$.

\bigskip

For notational convenience we assume without loss of generality that
$t\geq2$. We will repeatedly use in the sequel that, for $a>0$ and
$b\in\mathbb{R}$, there exists a constant $C>0$ such that 
\begin{equation} \label{repeatedly}
\int_2^t s^a(\log s)^b\,ds\leq C t^{a+1}(\log t)^b.
\end{equation}

It suffices to prove that there exists a sufficiently large constant
$C_1$ such that if $|x|\geq C_1t^{\frac12} (\log t)^{\frac{1-\theta}4}$, then 
\[\frac{x}{|x|}\cdot u(x,t)\leq C_1 t^{\frac12}(\log t)^{-\frac\theta2}|x|^{-2}.\] This will imply that 
a fluid particle cannot escape the region $|x|\leq \mzece t^{\frac12} (\log t)^{\frac{1-\theta}4}$ for a 
sufficiently large constant $\mzece$. To see this we reason by contradiction, following the trajectory of a 
fluid particle in a region like $|x|\geq C_1t^{\frac12} (\log t)^{\frac{1-\theta}4}$; we note that 
$\frac{x}{|x|}\cdot u$ is the radial velocity of said particle and we integrate the resulting ODE using 
\eqref{repeatedly} to deduce that the trajectory will remain in the region $|x|\leq C_2t^{\frac12} 
(\log t)^{-\frac{\theta}6}$ for some $C_2>0$. The conclusion follows since $C_2t^{\frac12} 
(\log t)^{-\frac{\theta}6} \leq \mzece t^{\frac12} (\log t)^{\frac{1-\theta}4}$ for some $M_7>0$.

Let $\widetilde\omega(x,t)=\omega(x,t)\bigl(\int_\dom
\domtx\bigr)^{-1}$. Since the function $[1,\infty)\ni s\mapsto \sigma(s)=
s^2\log^\theta s$ is convex and the integral of vorticity is a conserved quantity, we deduce from the Jensen inequality that
\begin{equation*}
  \sigma \Bigl(\int_\dom |T(x)| \widetilde\omega(x,t)\, dx\Bigr)\leq \int_\dom |T(x)|^2\log^\theta|T(x)| \widetilde\omega(x,t)\, dx\leq \mcinci (1+t)\|\omega_0\|_{L^1}^{-1}.
\end{equation*}
Next, one easily checks that $\sigma^{-1}(\rho)\sim 2^{\frac\theta2}
\rho^{\frac12}(\log\rho)^{-\frac\theta2}$ as $\rho\to\infty$. Therefore,
there exists a constant $C_3$ such that
\begin{equation}\label{1stmoment}
\int_\dom |T(x)|\domtx\leq C_3t^{\frac12}(\log t)^{-\frac\theta2}.  
\end{equation}

We assume now that  $|x|\geq C_1  t^{\frac12} (\log t)^{\frac{1-\theta}4}$ for a sufficiently large constant $C_1 $ to be determined later. We decompose $u=v+\al H_\dom$ where
\begin{align*}
  v(x,t)
&=\int_\dom \left(\frac{[(T(x)-T(y))DT(x)]^{\perp}}{2\pi|T(x)-T(y)|^2}
-
\frac{[(T(x)-T(y)^{\ast})DT(x)]^{\perp}}{2\pi|T(x)-T(y)^{\ast}|^2}\right)
\omega(y,t) dy\\
&\equiv \int_\dom K(x,y) \omega(y,t) dy,
\end{align*}
where $K(x,y)$ is the kernel from the first integral in the formula
above. We first estimate $v$. From the identity
$\left|\frac{a}{|a|^2}-\frac{b}{|b|^2}\right|=\frac{|a-b|}{|a||b|}$,
we deduce that
\begin{equation*}
  |K(x,y)|\leq \|DT\|_{L^\infty}\frac{|T(y)-T(y)^{\ast}|}{2\pi
    |T(x)-T(y)||T(x)-T(y)^\ast|}
\leq \|DT\|_{L^\infty}\frac{2|T(y)|}{\pi
    |T(x)-T(y)||T(x)|}.
\end{equation*}
We used above that the constant $C_1 $ may be chosen sufficiently large so as to have 
$|T(x)|>2>2|T(y)^{\ast}|$. We infer that 
\begin{align*}
  |v(x,t)|
&\leq \frac{C}{|T(x)|}\int_\dom \frac{|T(y)|}{|T(x)-T(y)|}
\omega(y,t) dy\\
&= \frac{C}{|T(x)|}\int_{|T(x)-T(y)|>|T(x)|/2} \frac{|T(y)|}{|T(x)-T(y)|}
\omega(y,t) dy\\
&\hskip 6cm + \frac{C}{|T(x)|}\int_{|T(x)-T(y)|<|T(x)|/2} \frac{|T(y)|}{|T(x)-T(y)|}
\omega(y,t) dy\\
&\equiv I_1+I_2.
\end{align*}

Next, by \eqref{1stmoment} we have that 
\begin{equation*}
  I_1\leq \frac{C}{|T(x)|^2}\int_\dom |T(y)|
\omega(y,t) dy
\leq Ct^{\frac12}(\log t)^{-\frac\theta2}|T(x)|^{-2}
\end{equation*}
and from Lemma \ref{velbd} we deduce that 
\begin{equation*}
  I_2 \leq C
\Bigl(\int_{|T(y)|>|T(x)|/2} 
\omega(y,t) dy\Bigr)^{\frac12}.
\end{equation*}

Next we estimate the radial component of $\alpha H_{\Omega^c}$. 
From the explicit formula for the harmonic vector field $H_\dom$ given
in relation \eqref{harmform} and from 
Lemma \ref{confmap}, we infer that 
$$\frac{x}{|x|}\cdot H_\dom(x)=\frac{x}{|x|}\cdot\Bigl(\frac{x^\perp}{2\pi|x|^2}+O\bigl(|x|^{-2}\bigr)\Bigr)=O\bigl(|x|^{-2}\bigr).$$
According to Proposition \ref{prop2} below, if the constant $C_1$ is
large enough then the term $I_2$ is also
$O\bigl(|x|^{-2}\bigr)$, so we finally deduce that 
\begin{equation*}
 \frac{x}{|x|}\cdot u(x,t)\leq |v(x,t)|+ \al \frac{x}{|x|}\cdot
 H_\dom(x)
\leq Ct^{\frac12}(\log t)^{-\frac\theta2}|x|^{-2},
\end{equation*}
with $C=C_1$.

This completes the proof of Theorem \ref{th1} once Proposition \ref{prop2} below is proved.
\end{proof}

\begin{proposition}\label{prop2}
Under the hypothesis of Theorem \ref{th1} we have that, for all $k>0$ there exists a constant $\msapte >0$ such that
  \begin{equation*}
    \int_{|T(x)|>r}\domtx\leq \frac{\msapte }{r^k}
  \end{equation*}
for all $r>\msapte  t^{\frac12} (\log t)^{\frac{1-\theta}4}$.
\end{proposition}
\begin{proof}
As before we assume $t\geq 2$.

  Let us introduce the following approximation for the mass of
  vorticity in the region of interest:
  \begin{equation*}
    f_r(t)=\int_\dom\eta\Bigl(\frac{|T(x)|^2-r^2}{\la r^2}\Bigr)\domtx,
  \end{equation*}
where $\la=\la(r)\in(0,1)$ is to be chosen later and
$\eta(s)=\frac{e^s}{1+e^s}$ verifies that
\begin{equation}\label{eta}
  |\eta'(s)|\leq \min\{\eta(s),e^{-|s|}\}, \qquad |\eta''(s)|\leq \eta(s).
\end{equation}
To simplify notation, we set
$\eta_r(x)=\eta\Bigl(\frac{|T(x)|^2-r^2}{\la r^2}\Bigr)$ and
$\eta'_r(x)=\eta'\Bigl(\frac{|T(x)|^2-r^2}{\la r^2}\Bigr)$. Clearly,
it suffices to prove that $f_r(t)\leq r^{-k}$ for $r>\msapte
t^{\frac12} (\log t)^{\frac{1-\theta}4}$ and $\msapte$ sufficiently large.

We differentiate $f_r$ to obtain
\begin{multline*}
  f'_r(t)=\int_\dom \eta_r(x)\partial_t\omega
=-\int_\dom \eta_r(x)u\cdot\nabla \omega=\int_\dom u\cdot\nabla
\eta_r\omega\\
=\iint_{\dom\times\dom} \nabla\eta_r(x)\cdot K(x,y) \omega(x,t)\omega(y,t)\,dx\,dy.
\end{multline*}
We used above that $\nabla\eta_r\cdot H_\dom=0$ pointwise; this
follows after noticing that the first term is proportional to
$\nabla|T|$ while the second one is proportional to
$\nabla^\perp|T|$. Using the explicit formula for $K(x,y)$ we finally
obtain that
\begin{equation*}
  \pi f'_r(t)=\frac1{\la r^2}\iint_{\dom\times\dom}
  \eta'_r(x)[T(x)DT(x)]\cdot \bigl[\widetilde{K}(x,y)DT(x)\bigr]^\perp
  \omega(x,t)\omega(y,t)\,dx\,dy, 
\end{equation*}
where
\begin{equation*}
\widetilde{K}(x,y)=\frac{T(x)-T(y)}{|T(x)-T(y)|^2}-\frac{T(x)-T(y)^\ast}{|T(x)-T(y)^\ast|^2}.
\end{equation*}
From Lemma \ref{confmap}, one can deduce that $DT=\beta
Id+O\Bigl(\frac1{|T(x)|^2}\Bigr)$. On the other hand, 
Lemma \ref{velbd} implies that 
$
  \int_\dom |\widetilde{K}(x,y)|\omega(y,t)\,dy\leq C.
$
We infer that
\begin{equation*}
 \pi f'_r(t)-  \frac{\beta^2}{\la r^2}\iint_{\dom\times\dom}
  \eta'_r(x)T(x)\cdot \widetilde{K}(x,y)^\perp
  \omega(x,t)\omega(y,t)\,dx\,dy
\leq \frac{C}{\la r^2}\int_\dom \frac {\eta'_r(x)}{|T(x)|}\domtx
\end{equation*}
Using \eqref{eta}, for $|T(x)|<r/\sqrt2$ we bound $\frac {\eta'_r(x)}{|T(x)|}\leq
e^{-\frac1{2\la}}$ while for $|T(x)|\geq r/\sqrt2$ we bound   $\frac
{\eta'_r(x)}{|T(x)|}\leq  \eta_r(x)\sqrt2 /r$. This implies that
\begin{equation}\label{a1}
 \pi f'_r(t)\leq J+ \frac{C}{\la r^3}f_r(t)+\frac{C}{\la r^2}e^{-\frac1{2\la}},
\end{equation}
where
\begin{equation}
  \label{a2}
  \begin{split}
  J
&=\frac{\beta^2}{\la r^2}\iint_{\dom\times\dom}
  \eta'_r(x) \frac{T(x)\cdot[T(x)-T(y)]^\perp}{|T(x)-T(y)|^2}
  \omega(x,t)\omega(y,t)\,dx\,dy\\
&\hskip 2cm -\frac{\beta^2}{\la r^2}\iint_{\dom\times\dom}
  \eta'_r(x) \frac{T(x)\cdot[T(x)-T(y)^\ast]^\perp}{|T(x)-T(y)^\ast|^2}
  \omega(x,t)\omega(y,t)\,dx\,dy\\
&\equiv J_1+J_2.    
  \end{split}
\end{equation}

We bound first $J_2$. Clearly we may assume, without loss of generality, that $r^2>8$. 
We decompose the integral in $J_2$ into an integral over $\{|T(x)| \leq 2 \}$ and an integral over 
$\{|T(x)|\geq 2\}$ and we write $J_2=J_{21}+J_{22}$ for each portion. 

Now, if $|T(x)|\leq2 < r/\sqrt{2}$ we estimate
\begin{equation*}
\Bigl|\eta'_r(x)
\frac{T(x)\cdot[T(x)-T(y)^\ast]^\perp}{|T(x)-T(y)^\ast|^2}\Bigr|
=\Bigl|\eta'_r(x) \frac{T(y)^*\cdot[T(x)-T(y)^\ast]^\perp}{|T(x)-T(y)^\ast|^2}\Bigr|
\leq \frac{e^{-\frac1{2\la}}}{|T(x)-T(y)^\ast|}  
\end{equation*}
Hence, by Lemma \ref{velbd}, we have 
\[|J_{21}|\leq C\frac{\beta^2}{\lambda r^2}e^{-\frac{1}{2\lambda}}.\]

Next, if $|T(x)|>2$ then, since $|T(y)^\ast|<1$, we find 
\begin{equation*}
 \frac{\bigl|T(x)\cdot[T(x)-T(y)^\ast]^\perp\bigr|}{|T(x)-T(y)^\ast|^2}
=\frac{\bigl|T(x)\cdot[T(y)^\ast]^\perp\bigr|}{|T(x)-T(y)^\ast|^2}
\leq \frac4{|T(x)|}.
\end{equation*}
Hence, after further decomposing the integral in $J_{22}$, we repeat the argument used to obtain
\eqref{a1} to find
\[|J_{22}| \leq  \frac{C}{\lambda r^2}e^{-\frac{1}{2\lambda}} + \frac{C}{\lambda r^3}f_r(t).\]

Putting these two estimates together we infer that 
\begin{equation}\label{a3}
|J_2|\leq \frac{C}{\la r^3}f_r(t)+\frac{C}{\la r^2}e^{-\frac1{2\la}}.
\end{equation}

We now turn to the estimate of $J_1$.
Using the change of variables $(x,y)\leftrightarrow (y,x)$ we next write
\begin{equation}\label{a4}
  J_1= \frac{\beta^2}{2\la r^2}\iint_{\dom\times\dom}
  L(x,y)
  \omega(x,t)\omega(y,t)\,dx\,dy
\end{equation}
where
\begin{equation*}
  L(x,y)=[\eta'_r(x)-\eta'_r(y)] \frac{T(x)\cdot[T(x)-T(y)]^\perp}{|T(x)-T(y)|^2}
\end{equation*}

We divide the domain of integration $\dom\times\dom$ in three
(overlapping) pieces:
$A=\{\bigl||T(x)|-|T(y)|\bigr|\geq r/4\}$,
$B=\{|T(x)|\text{ and }|T(y)|\not\in(r/2,3r/2)\}$ and 
$C=\{|T(x)|\text{ and }|T(y)|\in[r/4,5r/4]\}$. Clearly  $\dom\times\dom=A\cup B\cup C$.

For $(x,y)\in A$ we use \eqref{loops1} and \eqref{eta} to bound
\begin{equation*}
  |L(x,y)|\leq [\eta'_r(x)+\eta'_r(y)]  \min[|T(x)|,|T(y)|] \frac4r
\leq [\eta_r(x)|T(y)|+\eta_r(y)|T(x)|]\frac4r,
\end{equation*}
so that, in view of \eqref{1stmoment} 
\begin{equation}\label{a5}
 \frac{\beta^2}{2\la r^2}\iint_{A}
  |L(x,y)|
  \omega(x,t)\omega(y,t)\,dx\,dy 
\leq C\frac{f_r}{\la r^3}t^{\frac12}(\log t)^{-\frac\theta2}. 
\end{equation}

In the region $B$, as a consequence of \eqref{eta} both $\eta'_r(x)$
and $\eta'_r(y)$ are less than $e^{-\frac3{4\la}}$. Therefore, by
Lemma \ref{velbd}
\begin{equation}
  \label{a6}
  \begin{split}
 \frac{\beta^2}{2\la r^2}\iint_{B}
  |L(x,y)|
&
  \omega(x,t)\omega(y,t)\,dx\,dy \\
&\leq C\frac{e^{-\frac3{4\la}}}{\la r^2} \int_\dom |T(x)| \Bigl(\int_\dom\frac1{|T(x)-T(y)|}\omega(y,t)\,dy\Bigr)\domtx\\
&\leq C\frac{e^{-\frac3{4\la}}}{\la r^2} \int_\dom |T(x)| \domtx\\
&\leq C\frac{e^{-\frac3{4\la}}}{\la r^2} t^{\frac12}(\log t)^{-\frac\theta2}. 
  \end{split}
\end{equation}

Next, by the mean value theorem, relation \eqref{eta}  and since
 $\eta$ is increasing, one has that
$$|\eta'_r(x)-\eta'_r(y)|=\frac{\bigl||T(x)|^2-|T(y)|^2\bigr|}{\la
  r^2}|\eta''(\xi)|\leq\frac{\bigl||T(x)|^2-|T(y)|^2\bigr|}{\la
  r^2}[\eta_r(x)+\eta_r(y)],$$ 
for some $\xi$ between
$\frac{|T(x)|^2-r^2}{\la r^2}$ and $\frac{|T(y)|^2-r^2}{\la
  r^2}$. Therefore, in view of \eqref{loops1}, for $(x,y)\in C$ we can bound
\begin{equation*}
  |L(x,y)|\leq \frac{\bigl| |T(x)|^2-|T(y)|^2\bigr|}{\la
  r^2}[\eta_r(x)+\eta_r(y)] \frac{\min[|T(x)|,|T(y)|]}{ |T(x)-T(y)|}
\leq \frac{25}{8\la}[\eta_r(x)+\eta_r(y)] 
\end{equation*}
so that 
\begin{equation}\label{a7}
 \frac{\beta^2}{2\la r^2}
\iint_{C}
 | L(x,y)|
  \omega(x,t)\omega(y,t)\,dx\,dy 
\leq \frac{{25}\beta^2}{8\la^2 r^2}f_r(t)\int_{|T(x)|>\frac r4}\omega(x,t)\,dx
\leq \frac {Ctf_r}{\la^2 r^4\log^\theta r}
\end{equation}
where we used \eqref{momint} to bound the mass of vorticity in the
region $|T(x)|>r/4$ by $\mcinci (1+t)16 r^{-2}\bigl[\log(r/4)\bigr]^{-\theta}$.

Collecting relations \eqref{a1} to \eqref{a7} we finally deduce the following differential
inequality verified by $f_r$:
\begin{equation*}
  f'_r(t)\leq Cf_r\Bigl(\frac{t^{\frac12}}{\la r^3 (\log
      t)^{\frac\theta2}} +\frac t{\la^2 r^4\log^\theta  
      r}\Bigr)+C\frac{e^{-\frac1{2\la}}t^{\frac12}}{\la r (\log
      t)^{\frac\theta2}},\qquad t\geq2.
\end{equation*}

Since $\omega(\cdot,2)$ is compactly supported, choosing $r$ sufficiently large one has that $f_r(2)\leq
e^{-\frac1{2\la}}\|\omega_0\|_{L^1}$. Applying the Gronwall lemma in
the inequality above for $t\geq2$ and using \eqref{repeatedly}, we deduce that there exists a
constant $C_1$ such that
\begin{equation*}
f_r(t)\leq C_1\Bigl(1+ \frac{t^{\frac32}}{\la r (\log
  t)^{\frac\theta2}} \Bigr)\exp\Bigl( \frac{C_1t^{\frac32}}{\la r^3 (\log
      t)^{\frac\theta2}} +\frac {C_1t^2}{\la^2 r^4\log^\theta
      r}  -\frac1{2\la}\Bigr)
\end{equation*}

We finally choose $\la =\frac1{4n\log r}$ and $r$ sufficiently large such that
\begin{equation*}
 \frac{C_1t^{\frac32}}{\la r^3 (\log
      t)^{\frac\theta2}} +\frac {C_1t^2}{\la^2 r^4\log^\theta
      r} \leq \frac1{4\la}. 
\end{equation*}
It is a straightforward calculation to check that the above
restriction is implied by the condition $r>C_2 t^{\frac12} 
(\log t)^{\frac{1-\theta}4}$ for some large constant $C_2$. For this choice one
has that
\begin{equation*}
  \exp\Bigl( \frac{Ct^{\frac32}}{\la r^3 (\log
      t)^{\frac\theta2}} +\frac {Ct^2}{\la^2 r^4\log^\theta
      r}  -\frac1{2\la}\Bigr) \leq r^{-n}
\end{equation*}
For a large enough constant $n$ this term dominates $C_1+ \frac{C_1t^{\frac32}}{\la r (\log
  t)^{\frac\theta2}} $ and we obtain $f_r(t)\leq r^{-k}$ as desired.
\end{proof}

\section{Even vorticity on the exterior of the disk}

We consider now the case when $\Omega=D(0,1)$. The conformal map $T$
is simply the identity, so that the Biot-Savart law becomes:
\begin{equation}
  \label{0}
2\pi u(x,t)=\int\Bigl[\frac{(x-y)^\perp}{|x-y|^2}-
\frac{(x-y^\ast)^\perp}{|x-y^\ast|^2}\Bigl]\omega(y,t)\,dy
+\alpha\frac{x^\perp}{|x|^2},  
\end{equation}
where $\alpha$ is a real constant and $y^\ast=\frac{y}{|y|^2}$. It is
an easy calculation to show that the moment of inertia is conserved.
We  prove the following theorem:
\begin{theorem}\label{t1}
Suppose that the initial vorticity is nonnegative and even 
($\omega_0(-x)=\omega_0(x)\ \forall x$). There exists a constant $\mopt >0$
such that
$\supp \omega(\cdot,t)\subset
\bigl\{|x|\leq \mopt \bigl[(1+t)\log (2+t)\bigr]^{1/4}\bigr \}$ for some constant $\mopt $.
\end{theorem}
\begin{proof}
The proof follows the lines of the argument given in \cite{ISG99} with
the modifications induced by the new terms in the Biot-Savart
law. Indeed, the two keys facts used in \cite{ISG99} are also true in
this situation: both the moment of inertia and the center of mass are
conserved quantities (the center of mass is at 0 since the vorticity is even). It
suffices to assume that $t\geq2$ and to prove that there exists a constant $C_1$ such that 
\begin{equation}
  \label{1}
\frac{x}{|x|}\cdot u(x,t)\leq \frac{C_1}{|x|^3},  
\end{equation}
for any $|x|\geq C_1 (t\log t)^{1/4}$ such that $x$ is a point in the support of 
$\omega(\cdot,t)$ whose distance to the origin is the largest possible.

Note that $\omega(-x,t)=\omega(x,t)$ for all $x,t$. Next, 
let $x$ be a point as described above. Then,
\begin{multline*}
2\pi \frac{x}{|x|}\cdot u(x,t)= \frac{x}{|x|}\cdot  
\int\Bigl(\frac{{y^\ast}^\perp}{|x-{y^\ast}|^2}-
\frac{y^\perp}{|x-y|^2}\Bigl)\omega(y,t)\,dy
=\frac{x}{|x|^3}\cdot  
\int\bigl({y^\ast}^\perp-y^\perp\bigl)\omega(y,t)\,dy\\
+\frac{x}{|x|}\cdot  
\int\biggl[{y^\ast}^\perp\Bigl(\frac{1}{|x-{y^\ast}|^2}-\frac{1}{|x|^2}\Bigl)
-y^\perp\Bigl(\frac{1}{|x-y|^2}-\frac{1}{|x|^2}\Bigl)\biggl]\omega(y,t)\,dy.
\end{multline*}
For symmetry reasons,
\begin{equation}
  \label{2}
\int\bigl({y^\ast}^\perp-y^\perp\bigl)\omega(y,t)\,dy=0.  
\end{equation}
We also observe that, for any $x$ and $z$ we have 
\begin{equation*}
\frac{1}{|x-z|^2}-\frac{1}{|x|^2}=\frac{z\cdot(2x-z)}{|x|^2|x-z|^2}.  
\end{equation*}
We deduce that
\begin{equation}\label{finalI12}
2\pi \frac{x}{|x|}\cdot u(x,t)=  
\underset{I_1}{\underbrace{\int x\cdot{y^\ast}^\perp
\frac{{y^\ast}\cdot(2x-{y^\ast})}
{|x|^3|x-{y^\ast}|^2}\omega(y,t)\,dy}}
-\underset{I_2}{\underbrace{\int x\cdot y^\perp\frac{y\cdot(2x-y)}
{|x|^3|x-y|^2}\omega(y,t)\,dy}}.
\end{equation}
Recall that, by assumption, we always have that
$|y|\leq|x|$ in the integrands above. To bound the first integral, note that $|{y^\ast}|<1$ so that 
$|x-{y^\ast}|\geq|x|/2$ and $|2x-{y^\ast}|\leq 3|x|$. Consequently,
\begin{equation}\label{finalI1}
I_1\leq \frac{12}{|x|^3}\int |{y^\ast}|^2\omega(y,t)\,dy
\leq \frac{12}{|x|^3}\int \omega(y,t)\,dy.  
\end{equation}
We now go to the estimate of $I_2$ and split the domain of integration
in two pieces: $\{|y|<|x|/2\}$ and $\{|x|/2\leq|y|\leq |x|\}$. For
$|y|<|x|/2$, we still have that $|x-{y}|\geq|x|/2$ and $|2x-{y}|\leq
3|x|$, so as in the estimate for $I_1$, we can bound
\begin{equation*}
  \Bigl|x\cdot y^\perp\frac{y\cdot(2x-y)}
{|x|^3|x-y|^2} \Bigr|\leq \frac{12|y|^2}{|x|^3}\qquad \forall |y|<\frac{|x|}2.
\end{equation*}
Next, for $|x|/2\leq|y|\leq |x|$ we write
\begin{equation*}
 \Bigl|x\cdot y^\perp\frac{y\cdot(2x-y)}
{|x|^3|x-y|^2} \Bigr|
= \Bigl|(x-y)\cdot y^\perp\frac{y\cdot(2x-y)}
{|x|^3|x-y|^2} \Bigr|  
\leq \frac{|y|^2|2x-y|}{|x|^3|x-y|}\leq \frac3{|x-y|}.
\end{equation*}
We infer that
\begin{equation}
  \label{finalI2}  
\begin{split}
  |I_2|
&\leq \frac{12}{|x|^3}\int
  |y|^2\omega(y,t)\,dy+3\int_{|y|\geq|x|/2} \frac1{|x-y|}\omega(y,t)\,dy\\
&\leq \frac{12}{|x|^3}\int
  |y|^2\omega_0(y)\,dy+C \|\omega\|_{L^\infty}^{\frac12}\Bigl(\int_{|y|\geq|x|/2} \omega(y,t)\,dy\Bigr)^{\frac12},
\end{split}
\end{equation}
where we used that the moment of inertia is conserved and Lemma \ref{velbd}. 

Clearly relation \eqref{1} now follows from \eqref{finalI12},
\eqref{finalI1}, \eqref{finalI2} and from Lemma \ref{l1} below.
This completes the proof of Theorem \ref{t1} once Lemma \ref{l1} 
is proved.
\end{proof}
\begin{lemma}\label{l1}
For every $k\in\mathbb{N}$ there exists a constant $\mnoua(k)$ such that
\begin{equation*}
\int_{|x|\geq r}\omega(x,t)\,dx\leq\frac{\mnoua }{r^k}\quad\text{for}\quad
r\geq \mnoua(t\log t)^{1/4}.  
\end{equation*}
\end{lemma}
\begin{proof}
Let
\begin{equation*}
  f_r(t)=\int\varphi_r(x)\omega(x,t)\,dx,
\end{equation*}
where 
\begin{equation*}
\varphi_r(x)=\eta\Bigl(\frac{|x|^2-r^2}{\lambda r^2}\Bigr),
\quad \eta(s)=\frac{e^s}{1+e^s},  
\end{equation*}
and $\lambda=\lambda(r)\leq 1$ is a positive function to be chosen later. It suffices to prove that $f_r(t)\leq\frac{C(k)}{r^k}$.

As in the proof of Proposition \ref{prop2}, we deduce that 
$f'_r(t)$ can be written under the form
\begin{equation}\label{2a}
\begin{split}
2\pi f'_r(t)&=\iint\Bigl(\frac{{y^\ast}^\perp}{|x-{y^\ast}|^2}
-\frac{y^\perp}{|x-y|^2}\Bigr) \cdot\nabla\varphi_r(x)
\omega(x,t)\omega(y,t)\,dx\,dy\\  
&=\underset{J_1}{\underbrace{\iint\Bigl(\frac{1}{|x-{y^\ast}|^2}
-\frac{1}{|x|^2}\Bigr){y^\ast}^\perp \cdot\nabla\varphi_r(x)
\omega(x,t)\omega(y,t)\,dx\,dy}}\\
&\qquad\qquad\qquad-\underset{J_2}{\underbrace{\iint\Bigl(\frac{1}{|x-y|^2}
-\frac{1}{|x|^2}\Bigr)y^\perp \cdot\nabla\varphi_r(x)
\omega(x,t)\omega(y,t)\,dx\,dy}}
\end{split}
\end{equation}
where we have used relation \eqref{2}.

If we follow the analysis from relation (8) to the bottom of page 1720 of \cite{ISG99}, then we find that the second term can be estimated by
\begin{equation}\label{2b}
J_2\leq C\frac{f_r}{\lambda^2 r^4}+\frac{C}{\lambda r^2}
e^{-\frac{1}{2\lambda}}.  
\end{equation}
It remains to estimate 
\begin{multline*}
J_1=\iint_{|x|\leq r/2}\underset{L_1(x,y)}{\underbrace{
\Bigl(\frac{1}{|x-{y^\ast}|^2}
-\frac{1}{|x|^2}\Bigr)({y^\ast}-x)^\perp \cdot\nabla\varphi_r(x)}}
\omega(x,t)\omega(y,t)\,dx\,dy \\
+\iint_{|x|\geq r/2}\underset{L_2(x,y)}{\underbrace{
\frac{{y^\ast}\cdot(2x-{y^\ast})}{|x|^2|x-{y^\ast}|^2}
{y^\ast}^\perp \cdot\nabla\varphi_r(x)}}
\omega(x,t)\omega(y,t)\,dx\,dy. 
\end{multline*}
From the definition of $\varphi_r$, we have that
\begin{equation}
  \label{3}
\nabla\varphi_r(x)=\frac{2x}{\lambda r^2}\eta'\Bigl(
\frac{|x|^2-r^2}{\lambda r^2}\Bigr).  
\end{equation}
For $|x|\leq r/2$, one has that $\frac{|x|^2-r^2}{\lambda r^2}\leq 
-\frac{1}{2\lambda}$. Since $|\eta'(s)|\leq e^s$, it follows that 
\begin{equation*}
\Bigl| \eta'\Bigl(
\frac{|x|^2-r^2}{\lambda r^2}\Bigr)\Bigr|\leq e^{-\frac{1}{2\lambda}} 
\end{equation*}
so that, using also that $|x|, r\geq 1$,
\begin{equation*}
|L_1(x,y)|\leq \Bigl(\frac{1}{|x-{y^\ast}|^2}
+\frac{1}{|x|^2}\Bigr) |x-{y^\ast}|\frac{2|x|}{\lambda r^2}
e^{-\frac{1}{2\lambda}}
\leq \frac{C e^{-\frac{1}{2\lambda}}}{\lambda r|x-{y^\ast}|}
+\frac{C}{\lambda} e^{-\frac{1}{2\lambda}}   
\end{equation*}
which, in view of Lemma \ref{velbd}, implies that
\begin{equation}
  \label{4}
\Bigl|\iint_{|x|\leq r/2}L_1(x,y)
\omega(x,t)\omega(y,t)\,dx\,dy\Bigr|
\leq  \frac{C}{\lambda} e^{-\frac{1}{2\lambda}}.     
\end{equation}
To estimate $L_2(x,y)$, we use again equation \eqref{3} and the facts
that $|{y^\ast}|<1$ and $|\eta'|\leq \eta$ to deduce that, 
for $|x|\geq r/2$,
\begin{align*}
|L_2(x,y)|&=\Bigl|\frac{{y^\ast}\cdot(2x-{y^\ast})}{|x|^2|x-{y^\ast}|^2}
{y^\ast}^\perp \cdot\frac{2x}{\lambda r^2}\eta'\Bigl(
\frac{|x|^2-r^2}{\lambda r^2}\Bigr)\Bigr| \\ 
&\leq \frac{C}{|x|^2\lambda r^2}\eta\Bigl(
\frac{|x|^2-r^2}{\lambda r^2}\Bigr)\\
&\leq \frac{C}{\lambda r^4}\varphi_r(x).
\end{align*}
We infer that
\begin{equation}
  \label{5}
\Bigl|\iint_{|x|\geq r/2}L_2(x,y)
\omega(x,t)\omega(y,t)\,dx\,dy\Bigr|
\leq  \frac{C}{\lambda r^4}f_r(t).  
\end{equation}
Relations \eqref{4} and \eqref{5} now imply that
\begin{equation*}
|J_1|\leq \frac{C}{\lambda} e^{-\frac{1}{2\lambda}}
+ \frac{C}{\lambda r^4}f_r(t). 
\end{equation*}
Combining this with \eqref{2a} and \eqref{2b} we get
\begin{equation*}
f'(r)\leq \frac{Cf_r}{\lambda^2 r^4}+ 
\frac{C}{\lambda} e^{-\frac{1}{2\lambda}}. 
\end{equation*}
After integration, we obtain
\begin{equation*}
f_r(t)\leq (f_r(2)
+\lambda r^4  e^{-\frac{1}{2\lambda}})e^{\frac{Ct}{\lambda^2 r^4}}
\leq e^{\frac{Ct}{\lambda^2 r^4}-\frac{C_1}{\lambda}}
(C+\lambda r^4). 
\end{equation*}
Lemma \ref{l1} now follows by choosing $\lambda\log r$ sufficiently
small as we did at the end of the proof of Proposition \ref{prop2}.
\end{proof}

Let us conclude with the remark that, although we assumed that $\omega_0$ was smooth, our results hold, with minor modifications, if $\omega_0 \in L^p_c$ for some $p>2$.

{\footnotesize  {\it Acknowledgments:} 
This research has been supported in part by the UNICAMP Differential Equations PRONEX, FAPESP grant \# 00/02097-1,
FAEP grant \# 0285/01, CNPq grant \# 472.504/2004-5 and by the CNPq-CNRS Brazil-France Cooperation. The first author wishes to thank IMECC-UNICAMP for its hospitality, while the second and third authors are grateful to the Mathematics Departments of
Univ. de Rennes I and of Univ. de Lyon I for their hospitality. This work was completed during the Special Semester in Fluid Mechanics at the Centre Interfacultaire Bernoulli, EPFL; the authors wish to express their gratitude for the hospitality received.}

\noindent
{\sc
Drago\c{s} Iftimie\\
Institut Camille Jordan, Universit{\'e} Lyon 1,\\
B{\^a}t. Braconnier, 21 av. Claude Bernard, \\
69622 Villeurbanne cedex, France
\\}
{\it E-mail address:} dragos.iftimie@univ-lyon1.fr

\vspace{.1in}
\noindent
{\sc
Milton C. Lopes Filho\\
Departamento de Matematica, IMECC-UNICAMP.\\
Caixa Postal 6065, Campinas, SP 13083-970, Brasil
\\}
{\it E-mail address:} mlopes@ime.unicamp.br

\vspace{.1in}
\noindent
{\sc 
Helena J. Nussenzveig Lopes\\
Departamento de Matematica, IMECC-UNICAMP.\\
Caixa Postal 6065, Campinas, SP 13083-970, Brasil
\\}
{\it E-mail address:} hlopes@ime.unicamp.br

\end{document}